\documentclass{article} 
\usepackage[top=3.5cm,right=3cm,bottom=3.5cm,left=3cm]{geometry}
\usepackage[english]{babel} 
\usepackage{amssymb}
\usepackage{amsmath}
\usepackage{txfonts}
\usepackage{mathdots}
\usepackage[classicReIm]{kpfonts}
\usepackage{graphicx}

\begin{document}


\noindent 

\noindent 

\noindent 

\noindent 
\begin{center}
\noindent \textbf{\large Comparison of solutions resulted from direct problems formulated as FRE }
\end{center}
\bigskip
\noindent \textbf{}
\begin{center}

\noindent \textbf{Amin Ghodousian${}^{a,}$}${}^{}$\footnote{\ Corresponding\ author$  \\$Email\ address:\ \ $  $a.ghodousian@ut.ac.ir\ (Amin\ Ghodousian).}\textbf{, Sara Zal${}^{b}$}

\end{center}
\bigskip
\noindent {\footnotesize ${}^{a,}$Faculty of Engineering Science, College of Engineering, University of Tehran, P.O.Box 11365-4563, Tehran, Iran.\noindent\medskip }

\noindent {\footnotesize ${}^{b}$Department of Engineering Science, College of Engineering, University of Tehran, Tehran, Iran.\textbf{}}

\noindent \textbf{}

\noindent \textbf{Abstract}
\vskip 0.2in
\noindent In this paper, we investigate direct solution of FRE and compare their results with expected real consequences. We give an applied example formulated by a FRE problem and show that FRE defined by maximum t-conorm and an arbitrary t-norm can yield different interpretations for our example. A necessary condition and a sufficient condition are presented that guarantee FRE defined by maximum t-conorm and minimum t-norm attain the same solutions as human mind does. Also, we present a t-conorm and use it instead of maximum t-conorm in FRE to obtain solutions with highest similarity with real ones. Finally, we show that under some conditions, FRE defined by any t-conorm and any t-norm may find a solution which is not reasonable.


\noindent 
\vskip 0.2in
\noindent \textbf{1. Introduction}
\vskip 0.2in
\noindent Fuzzy relational equations (FRE) are a generalization of Boolean relation equations. At first, Sanchez [39] investigated and developed this theory and its applications. Generally, fuzzy set theory has a number of properties that make it suitable for formulizing the uncertain information upon which many applied concepts such as medical diagnosis and treatment are usually based. Because of this fact the theory of FRE was originally applied in problems of the medical diagnosis [39]. Later, The concept of FRE is used in many problems such as system analysis[38], decision making[1], fuzzy controller [10], fuzzy modeling[42], fuzzy analysis, especially fuzzy arithmetic[10]. In [21] it has been established that the majority of fuzzy inference systems can be implemented using FRE. Pedrycz [30] categorized and extended two ways of the generalizations of FRE in terms of sets under discussion and various operations which are taken into account. Since then , many theoretical improvements have been investigated and many applications have been presented. For instance, we can refer to [3,5,7,11-14,16,17,22,25,31-33,35,37,41,43]. Klement et al. [18-20] presented the basic analytical and algebraic properties of triangular norms and important classes of fuzzy operators generalization such as Archimedean, strict and nilpotent t-norms. Pedrycz and Vasilakos [34] converted a highly dimensional relational equation into a serious of single input FRE. In [36] author demonstrate how problems of interpolation and approximation of fuzzy functions are connected with solvability of systems of FRE. He focused on the problem of approximate solvability of a system of FRE. In [40], FRE has been extended to the setting of interval-valued FRE with a max-t-norm composition and three types of solution sets have been proposed. Markovskii showed that solving max-product FRE is closely related to the covering problem which is an NP-hard problem [24]. In [2] Chen and Wang designed an algorithm for obtaining the logical representation of all minimal solutions. They showed that a polynomial time algorithm to find all minimal solutions of FRE with max-min composition may not exist. 
\vskip 0.1in
\noindent   An interesting application of the fuzzy relations theory is in the field of image processing [9] and [6]. Also, the various types of FRE with continuous triangular norms were used for compression / decompression of images [4,15,28,23]. In [4] the image was divided in blocks and then the FRE with a t-norm were utilized for compression each block where results were obtained using the Lukasiewicz t-norm . In [26] any monochromatic image is interpreted as a fuzzy relation in which the entries are the normalized values of the pixels. This method is based essentially on the fact that the reconstructed images is obtained as the greatest or the smallest solution of a system of fuzzy equations. Also, authors showed processes for coding/decoding color images in the RGB and YUV spaces by using FRE of max-t type where t is the Yager t-norm. In [27] authors used particular FRE of compression / decompression of color images in the RGB and YUV spaces .Nobuhara et al. [29] formulated and solved a problem of image reconstruction using eigen fuzzy sets. They proposed two algorithms of generating eigen fuzzy sets used in the reconstruction process.Di Nola and Russo [8] focused on the algebraic structures of such these problems by considering the Lukasiewicz transform as a residuated map.
\vskip 0.1in
\noindent    In this paper, we firstly focus on a direct problem defined by the most familiar form of fuzzy relational equations which is formulated as follows:

\noindent 
\begin{equation} \label{GrindEQ__1_} 
xo \, A=b 
\end{equation}

\noindent Where$I\, =\, \{ \, 1\, ,\, 2\, ,\, ...\, ,\, n\, \} $, and$J\, =\, \{ \, 1\, ,\, 2\, ,\, ...\, ,\, m\, \} $.$A\, =\, (\, a_{ij} \, )_{n\times m} $ is a fuzzy matrix in which $0\, \, \le \, \, a_{ij} \, \, \le \, \, 1$,$\forall \, i\in \, I$ and $\forall j\in J\, $. $x$ is a given $n$-dimensional fuzzy vector, and $b$ is an unknown $m$-dimensional fuzzy vector in which $0\le x_{i} \le 1$ ,$\forall i\in I$, and $0\le b_{j} \le 1$, $\forall j\in J$, respectively. Also, "$o $" is an arbitrary t-norm, and $a_{j} $ is the $j$th column of matrix$A$. Equivalently, Problem(1) can be written with more details as follows:

\noindent                                                                                                                     
\begin{equation} \label{GrindEQ__2_} 
xo \, a_{j} \, =\, \mathop{\max }\limits_{i=1}^{n} \{ \, t(\, x_{i} ,a_{ij} )\, \} \, =\, b_{i}    ,\forall j\in J 
\end{equation}

\noindent   We present an applied example that is formulated by (2) as a direct problem in which $x$ and $A$ are given and $b$ is unknown. We present a necessary condition and a sufficient condition in which direct solution of Problem (2) defined by minimum t-norm leads to a solution coinciding with consequence resulted from human mind. Also, we show that direct Problem (2) defined by each t-norm cannot yield a rational vector$b$ for our example. In latter case, provided with some assumptions hold, we present a t-conorm (instead of maximum operator in (2)) by which we can solve direct Problem (2) (in which $t(a,b)=\min (a,b)$) with acceptable results and we find a vector $b$ that is exactly the same as rational consequences which are obtained by human mind. Also, we investigate conditions in which direct Problem (2) defined by each t-norm and each t-conorm cannot yield a reasonable solution.

\noindent 

\noindent   In section 2, we give an example to prove that direct Problem (2) may result in a vector $b$ which is too farther from the expected rational conclusions. In section 3, we present a t-conorm which is constructed by the convex combination of drastic sum and maximum operators. Then, we solve direct Problem (2) defined by minimum t-norm in which foregoing t-conorm takes the place of maximum operator, and we show that vector $b$obtained by the new problem has the highest similarity with the expected rational conclusions. Finally, conditions are presented in which any kind of Problem (2) (defined by arbitrary t-conorm and t-norm) does not result in a rational solution.

\noindent 
\vskip 0.2in
\noindent \textbf{2. An applied example: pea or melon?}
\vskip 0.2in
\noindent Consider a basket including several peas, strawberries, bananas, and melons. Also, consider a man who is aware of those four types of contents deciding to distinguish those objects only by touching them with closed eyes. By this assumption that the contents in the basket are just those four foregoing kinds, one will able to evaluate them in his mind by considering some properties(for instance, their weights). When one touches an object with sufficient large value of weight (for example he assigns 1 to the weight of that object from the interval [0,1]) he may argue that what is in his hands is melon. On the other hand, when one touches a melon, he can feel that what is in his hands has large value of weight among other contents included in the basket. Then, we can define a fuzzy relation in which, membership value of being melon and simultaneously having large value of weight, is too high. By similar reasons, we may make a table based on four objects; pea, strawberry, banana, melon, and three properties; long, Heavy, and voluminous as follows:   
\vskip 0.2in
\noindent 
\begin{center}
\begin{tabular}{|p{0.8in}|p{0.7in}|p{0.7in}|p{0.7in}|p{0.7in}|} \hline 
 &       Pea & Strawberry &     Banana &     Melon \\ \hline 
      Long & 0 & 0.3 & 0.98 & 0.7 \\ \hline 
     Heavy & 0.001 & 0.01 & 0.1 & 0.99 \\ \hline 
Voluminous & 0.004 & 0.042 & 0.3 & 1 \\ \hline 
\end{tabular}

\noindent\textbf{\footnotesize      Table 1. }{\footnotesize  Fuzzy relation between four objects; pea, strawberry, banana, melon, and three properties; long, heavy, voluminous.}
\end{center}
\noindent 

\noindent 

\noindent   By our assumption, an object selected from basket is melon iff its length = 0.7 or its weight = 0.99 or its volume = 1. Similar statements are true for being banana, being strawberry, and being pea. Furthermore, table 1 can be shown by a fuzzy relational matrix as follows:
\vskip 0.2in
\[\begin{array}{c} {} \\ {Length} \\ {Weight} \\ {Volume} \end{array}\begin{array}{l} {\, \, \, \, \, \, \, \, \, \, \, \, \, \, \, P\, \, \, \, \, \, \, \, \, \, \, \, S\, \, \, \, \, \, \, \, \, \, \, B\, \, \, \, \, \, \, \, \, \, M\, \, \, \, } \\ {\, \, \, \, \, \left[\, \begin{array}{cccc} {0} & {0.3} & {0.98} & {0.7} \\ {0.001} & {0.01} & {0.1} & {0.99} \\ {0.004} & {0.042} & {0.3} & {1} \end{array}\, \right]=A} \end{array}\] 
\vskip 0.1in
\noindent where$P,S,B,$and$M$denote pea, strawberry, banana, and melon, respectively.
\vskip 0.1in
\noindent   

\noindent    Now, suppose we select an object $x$ from the basket (with closed eyes) and we evaluate its properties as; $length(x)=x_{1} $, $weight(x)=x_{2} $, and $volume(x)=x_{3} $, where $x_{1} ,x_{2} ,x_{3} \in [0,1]$. We show these pieces of information by a 3-dimensional vector as shown below:
\vskip 0.2in
\[\begin{array}{l} {\, \, \, \, \, \, \, \, \, L\, \, \, \, \, \, \, \, W\, \, \, \, \, \, \, V} \\ {x=\left[\begin{array}{ccc} {x_{1} } & {x_{2} } & {x_{3} } \end{array}\right]} \end{array}\] 
\vskip 0.1in
\noindent where $L,W,$and $V$ denote length, weight, and volume, respectively.
\vskip 0.1in
\noindent 

\noindent    Then, by the composition of fuzzy relations, we can find the kind of this object by observing vector $b$ which is attained by solving direct Problem (2) as follows:
\vskip 0.2in
\[\begin{array}{l} {\, \, \, \, \, \, \, \, \, \, \, \, \, \, \, \, \, \, \, \, \, \, \, \, \, \, \, \, \, \, \, \, \, \, \, \, \, \, \, \, \, \, \, \, P\, \, \, \, \, \, \, \, \, \, \, \, S\, \, \, \, \, \, \, \, \, \, \, B\, \, \, \, \, \, \, \, \, \, M\, \, \, \, \, \, \, \, \, \, \, \, \, P\, \, \, \, \, \, S\, \, \, \, \, \, B\, \, \, \, \, \, \, M\, \, \, \, \, \, \, \, \, } \\ {\, \, \left[\begin{array}{ccc} {x_{1} } & {x_{2} } & {x_{3} } \end{array}\right]\, \, o \, \left[\, \begin{array}{cccc} {0} & {0.3} & {0.98} & {0.7} \\ {0.001} & {0.01} & {0.1} & {0.99} \\ {0.004} & {0.042} & {0.3} & {1} \end{array}\, \right]=\left[\begin{array}{cccc} {b_{1} } & {b_{2} } & {b_{3} } & {b_{4} } \end{array}\right]} \end{array}\] 
\vskip 0.1in
 Vector $b$determines the membership value of each contents in the basket to which object $x$ belongs. Precisely,$b_{1} $ is the degree of which $x$ is pea, $b_{2} $ is the degree of which $x$ is strawberry, $b_{3} $is the degree of which $x$ is banana, and $b_{4} $ is the degree of which $x$ is melon.

\noindent   In a certain experiment, we select an object $\dot{x}$ from the basket and find that$\dot{x}=\left[\begin{array}{ccc} {0.004} & {0.002} & {0.003} \end{array}\right]$. By comparing three properties of object$\dot{x}$ with that of pea, strawberry, banana, and melon, it is easy to see that $\dot{x}$ should be a pea. Therefore, we must have the same conclusion by finding vector$b$ through solving Direct problem (2), too. At first, we solve problem $\dot{x}o \, A=b$ in which "$o $" is minimum t-norm and
\vskip 0.1in
\[\begin{array}{l} {\, \, \, \, \, \, \, \, \, \, \, \, \, \, \, \, \, \, \, \, \, \, P\, \, \, \, \, \, \, \, \, \, \, \, S\, \, \, \, \, \, \, \, \, \, \, B\, \, \, \, \, \, \, \, \, \, M\, \, \, \, } \\ {\, \, \, \, A=\, \left[\, \begin{array}{cccc} {0} & {0.3} & {0.98} & {0.7} \\ {0.001} & {0.01} & {0.1} & {0.99} \\ {0.004} & {0.042} & {0.3} & {1} \end{array}\, \right]} \end{array}\] 
\vskip 0.1in
Easy calculations show that the result is the vector below:
\vskip 0.1in
\[\begin{array}{l} {\, \, \, \, \, \, \, \, \, \, \, \, \, \, \, \, P\, \, \, \, \, \, \, \, \, \, \, \, \, \, S\, \, \, \, \, \, \, \, \, \, \, \, \, B\, \, \, \, \, \, \, \, \, \, \, \, \, M} \\ {\, b=\, \, \left[\begin{array}{cccc} {0.003} & {0.004} & {0.004} & {0.004} \end{array}\right]} \end{array}\] 
\vskip 0.1in
  According to the vector $b$ above, $\dot{x}$ is a pea with membership value 0.003, and it is a melon with membership value 0.004 (?!). If we solve problem above with product t-norm, we will attain solutions as objectionable as previous one: 
\[\begin{array}{l} {\, \, \, \, \, \, \, \, \, \, \, \, \, \, \, \, \, \, \, P\, \, \, \, \, \, \, \, \, \, \, \, \, \, \, \, \, S\, \, \, \, \, \, \, \, \, \, \, \, \, \, \, \, \, \, B\, \, \, \, \, \, \, \, \, \, \, \, \, \, \, \, M} \\ {\, b=\, \, \left[\begin{array}{cccc} {0.000012} & {0.0012} & {0.00392} & {0.003} \end{array}\right]} \end{array}\] 
  The result will be too wonderful if we solve problem above with Lukasiewicz t-norm or drastic product. In both cases, we have:
\[\begin{array}{l} {\, \, \, \, \, \, \, \, \, \, \, P\, \, \, \, S\, \, \, \, B\, \, \, \, \, \, \, \, M} \\ {\, b=\, \, \left[\begin{array}{cccc} {0} & {0} & {0} & {0.003} \end{array}\right]} \end{array}\] 
  More wonderful result is attained when we try to distinguish an object $\ddot{x}=\left[\begin{array}{ccc} {0} & {0.001} & {0.004} \end{array}\right]$ via direct solution of Problem (2). By considering the components of table 1, it is readily attained that $\ddot{x}$ is exactly a pea. But, when Problem (2) is  defined by minimum, product, Lukasiewicz, and drastic product, yields vectors$b_{1} $, $b_{2} $, $b_{3} $, and $b_{4} $ below, respectively:
\[\begin{array}{l} {\, \, \, \, \, \, \, \, \, \, \, \, \, \, \, \, \, \, P\, \, \, \, \, \, \, \, \, \, \, \, \, \, S\, \, \, \, \, \, \, \, \, \, \, B\, \, \, \, \, \, \, \, \, \, \, \, \, M} \\ {\, b_{1} =\, \, \left[\begin{array}{cccc} {0.004} & {0.004} & {0.004} & {0.004} \end{array}\right]} \end{array}\] 

\[\begin{array}{l} {\, \, \, \, \, \, \, \, \, \, \, \, \, \, \, \, \, \, \, P\, \, \, \, \, \, \, \, \, \, \, \, \, \, \, \, \, \, \, \, \, S\, \, \, \, \, \, \, \, \, \, \, \, \, \, \, \, \, \, \, \, \, \, B\, \, \, \, \, \, \, \, \, \, \, \, \, \, \, \, \, \, \, \, M} \\ {\, b_{2} =\, \, \left[\begin{array}{cccc} {16\times 10^{-6} } & {168\times 10^{-6} } & {12\times 10^{-4} } & {4\times 10^{-3} } \end{array}\right]} \end{array}\] 

\[\begin{array}{l} {\, \, \, \, \, \, \, \, \, \, \, \, \, \, \, \, \, \, \, \, \, \, \, \, P\, \, \, \, S\, \, \, \, B\, \, \, \, \, \, M} \\ {\, b_{3} =b_{4} =\, \, \left[\begin{array}{cccc} {0} & {0} & {0} & {0.004} \end{array}\right]} \end{array}\]

\noindent    Actually, we expect $b_{1} \ge b_{2} \ge b_{3} \ge b_{4} $ by a rational justification in all cases above while results attained from Problem (2) have inverse ordering in all cases.
\vskip 0.1in
\noindent   In a general case, suppose $\stackrel{...}{x}=\left[\begin{array}{ccc} {\stackrel{...}{x}_{1} } & {\stackrel{...}{x}_{2} } & {\stackrel{...}{x}_{3} } \end{array}\right]$ is a selected object from the basket. Since $a_{i1} <a_{i4} $ for each$i\in \{ \, 1\, ,\, 2\, ,\, 3\, \} $, then by the non-decreasing property of t-norms we have $t(\stackrel{...}{x}_{i} \, ,\, a_{i1} )\le t(\stackrel{...}{x}_{i} \, ,\, a_{i4} )$ for each$i\in \{ \, 1\, ,\, 2\, ,\, 3\, \} $and for each t-norm.Therefore, $\mathop{\max }\limits_{i=1}^{3} \{ t(\stackrel{...}{x}_{i} \, ,\, a_{i1} )\} \le \mathop{\max }\limits_{i=1}^{3} \{ t(\stackrel{...}{x}_{i} \, ,\, a_{i4} )\} $which implies $b_{1} \le b_{4} $. Since latter inequality is satisfied for each $\stackrel{...}{x}\in [0,1]^{3} $ and for each t-norm, then problem (2) defined by each t-norm cannot truly distinguish between a pea and a melon.  

\noindent \textbf{}

\noindent \textbf{3. Modification of Problem (2)}
\vskip 0.1in
\noindent In this section, we investigate how human mind can solve example above with true solutions while Problem (2) does not. Also, we try to modify Problem (2) in a way which leads us to the same consequences as made by human mind. The process which is used in the mind for solving example above is based on the comparison, and comparison process itself is done by a difference operation. Actually, if we want to formulate the solution process of mind for our example with highest similarity with direct problem (2), we have:
\begin{equation} \label{GrindEQ__3_} 
1-\mathop{\min }\limits_{i=1}^{n} \{ \, \left|\, x_{i} -a_{ij} \, \right|\, \} =b_{j} \, \, \, \, ,\forall j\in J 
\end{equation} 
or equivalently;

\noindent 
\begin{equation} \label{GrindEQ__4_} 
\mathop{\max }\limits_{i=1}^{n} \{ \, 1-\left|\, x_{i} -a_{ij} \, \right|\, \} =b_{j} \, \, \, \, ,\forall j\in J 
\end{equation} 
in which $x_{i} $ and $a_{ij} $ are given$\forall i\in I$ and $\forall j\in J$, and $b_{j} $ is unknown $\forall j\in J$.

\noindent 
\vskip 0.1in
\noindent    If we solve example above for $\dot{x}=\left[\begin{array}{ccc} {0.004} & {0.002} & {0.003} \end{array}\right]$ and $\ddot{x}=\left[\begin{array}{ccc} {0} & {0.001} & {0.004} \end{array}\right]$ by direct problem \eqref{GrindEQ__3_} or \eqref{GrindEQ__4_}, we will have reasonable vectors $\dot{b}$, and $\ddot{b}$, respectively as follows:
\[\begin{array}{l} {\, \, \, \, \, \, \, \, \, \, \, \, \, \, \, \, \, P\, \, \, \, \, \, \, \, \, \, \, \, S\, \, \, \, \, \, \, \, \, \, \, \, \, \, B\, \, \, \, \, \, \, \, \, \, \, \, M} \\ {\, \dot{b}=\, \, \left[\begin{array}{cccc} {0.999} & {0.992} & {0.902} & {0.304} \end{array}\right]} \end{array}\] 
\[\begin{array}{l} {\, \, \, \, \, \, \, \, \, \, \, P\, \, \, \, \, \, \, \, \, S\, \, \, \, \, \, \, \, \, \, \, \, B\, \, \, \, \, \, \, \, \, \, M} \\ {\, \ddot{b}=\, \, \left[\begin{array}{cccc} {1} & {0.991} & {0.974} & {0.3} \end{array}\right]} \end{array}\]

\noindent    Here the question is whether there exists a t-norm by which vector $b$resulted from Problem (2) is the same as one attained from Problem \eqref{GrindEQ__4_}, or else by which t-norm, Problem (2) yields a vector $b$ with the most accurate approximation to vector $b$ in \eqref{GrindEQ__4_}. Lemma 1 and its corollary below show that not only vector $b$ resulted from (2) defined by an arbitrary t-norm is not generally equal to its corresponding vector in \eqref{GrindEQ__4_}, but also under some conditions, the amount of the accuracy for vector $b$ in (2) as an approximation solution for vector $b$ in \eqref{GrindEQ__4_} is bounded above. However, Problem (2) is not generally equal to Problem \eqref{GrindEQ__4_}, and also under some conditions, the accuracy of Problem (2) (defined by each t-norm) to approximate Problem \eqref{GrindEQ__4_} can not be desirably increased.

\noindent 
\vskip 0.1in
\noindent \textbf{Lemma 1.} Given Problem (2). Then,$\mathop{\max }\limits_{i=1}^{n} \{ \, t(\, x_{i} ,a_{ij} )\} \le $$\mathop{\max }\limits_{i=1}^{n} \{ \, \min (x_{i} ,a_{ij} )\, \} \le $$\mathop{\max }\limits_{i=1}^{n} \{ 1-\left|x_{i} -a_{ij} \right|\, \} $$,\forall j\in J$, where $t$ is an arbitrary t-norm. 
\vskip 0.1in
\noindent \textbf{Proof.} The first inequality is well known. For the second, fix a $j_{0} \in J$ and suppose$x_{i} \le a_{ij_{0} } $for some$i\in I$(otherwise, if$x_{i} >a_{ij_{0} } $, result is similarly attained). In this case, if$x_{i} \le 2a_{ij_{0} } -1$, then$x_{i} \le 1-(a_{ij_{0} } -x_{i} )\le a_{ij_{0} } $. Otherwise, if$x_{i} >2a_{ij_{0} } -1$, then$x_{i} \le a_{ij_{0} } <1-(a_{ij_{0} } -x_{i} )$. Therefore, in both cases, we have:
\begin{equation} \label{5} 
\, \min (x_{i} ,a_{ij_{0} } )\, \le 1-\left|x_{i} -a_{ij_{0} } \right|\,  
\end{equation} 
 Now, since last inequality holds for each$i\in I$, we have $\mathop{\max }\limits_{i=1}^{n} \{ \, \min (x_{i} ,a_{ij_{0} } )\, \} \le $$\mathop{\max }\limits_{i=1}^{n} \{ 1-\left|x_{i} -a_{ij_{0} } \right|\, \} $. Since $j_{0} \in J$ was arbitrary, then latter inequality is satisfied for each$j\in J$. 

\noindent 
\vskip 0.1in
\noindent \textbf{Corollary 1.} Given Problem (2)defined by minimum t-norm. Consider a fixed $j_{0} \in J$ and suppose $\min (x_{i} ,a_{ij_{0} } )>2\max (x_{i} ,a_{ij_{0} } )-1$, $\forall i\in I$. Then, $\mathop{\max }\limits_{i=1}^{n} \{ \, t(\, x_{i} ,a_{ij_{0} } )\} \le $ $\mathop{\max }\limits_{i=1}^{n} \{ \, \min (x_{i} ,a_{ij_{0} } )\, \} <$$\mathop{\max }\limits_{i=1}^{n} \{ 1-\left|x_{i} -a_{ij_{0} } \right|\, \} $.

\noindent 
\vskip 0.1in
\noindent   Corollary 1 shows that under some conditions, Problem (2) defined by each  t-norm may not even a good approximation of real solution resulted from Problem \eqref{GrindEQ__4_}. Lemma 2 part (a) below gives a sufficient condition by which Problem (2) defined by minimum t-norm finds rational vector $b$, exactly. Also, a necessary condition has been given in part (b) when Problem (2) defined by minimum t-norm provides a reasonable vector $b$ the same as we expect. Moreover, 

\noindent 
\vskip 0.1in
\noindent \textbf{Lemma 2.} Given Problem (2). 

\noindent a) If  $\max (x_{i} ,a_{ij} )=1$,$\forall i\in I$ and $\forall j\in J$, then $\mathop{\max }\limits_{i=1}^{n} \{ \, \min (x_{i} ,a_{ij} )\, \} =$ $\mathop{\max }\limits_{i=1}^{n} \{ 1-\left|x_{i} -a_{ij} \right|\, \} $,$\forall j\in J$ .

\noindent b) Suppose $\mathop{\max }\limits_{i=1}^{n} \{ \, \min (x_{i} ,a_{ij} )\, \} =$$\mathop{\max }\limits_{i=1}^{n} \{ 1-\left|x_{i} -a_{ij} \right|\, \} $,$\forall j\in J$. Then, for each $j\in J$ there exist at least one$i_{j} \in I$ such that $\max (x_{i_{j} } ,a_{i_{j} j} )=1$.
\vskip 0.1in
\noindent \textbf{Proof.} a) If $\max (x_{i} ,a_{ij} )=1$, then $\min (x_{i} ,a_{ij} )=$$1-\left|x_{i} -a_{ij} \right|\, $. Therefore, if $\max (x_{i} ,a_{ij} )=1$, $\forall i\in I$, we have $\min (x_{i} ,a_{ij} )=$$1-\left|x_{i} -a_{ij} \right|\, $, $\forall i\in I$ which implies part (a). \noindent (b) Consider a fixed $j_{0} \in J$ and suppose $\mathop{\max }\limits_{i=1}^{n} \{ \, \min (x_{i} ,a_{ij_{0} } )\, \} =$$\mathop{\max }\limits_{i=1}^{n} \{ 1-\left|x_{i} -a_{ij_{0} } \right|\, \} $.Let
\begin{center}
\noindent $\min (x_{i_{1} } ,a_{i_{1} j_{0} } )=\mathop{\max }\limits_{i=1}^{n} \{ \, \min (x_{i} ,a_{ij_{0} } )\, \} $and$1-\left|x_{i_{2} } -a_{i_{2} j_{0} } \right|$$=\mathop{\max }\limits_{i=1}^{n} \{ 1-\left|x_{i} -a_{ij_{0} } \right|\, \} $
\end{center}
\noindent  Then, by assumption of part (b), we have:
\[\min (x_{i_{1} } ,a_{i_{1} j_{0} } )=1-\left|x_{i_{2} } -a_{i_{2} j_{0} } \right|                                           \qquad(*)\] 
Now, by (5), we have:
\[\min (x_{i_{1} } ,a_{i_{1} j_{0} } )\le 1-\left|x_{i_{1} } -a_{i_{1} j_{0} } \right|\le \mathop{\max }\limits_{i=1}^{n} \{ 1-\left|x_{i} -a_{ij_{0} } \right|\, \} =1-\left|x_{i_{2} } -a_{i_{2} j_{0} } \right|\] 
Which together with (*) imply $\min (x_{i_{1} } ,a_{i_{1} j_{0} } )=$$1-\left|x_{i_{1} } -a_{i_{1} j_{0} } \right|$. Since $j_{0} \in J$ was arbitrary, last equality completes the proof.

\noindent 
\vskip 0.1in
\noindent \textbf{Lemma 3.} Consider Problem (2) and Suppose $i_{1} ,i_{2} \in I$ and $j\in J$.
\vskip 0.1in
\noindent  a) If we have $\min (x_{i_{1} } ,a_{i_{1} j} )\, =\, \min (x_{i_{2} } ,a_{i_{2} j} )=0$ or, if $\min (x_{i_{1} } ,a_{i_{1} j} )\, =0$ and$\, \min (x_{i_{2} } ,a_{i_{2} j} )\ne 0$,then:
\[S_{D} \{ \min (x_{i_{1} } ,a_{i_{1} j} )\, ,\, \min (x_{i_{2} } ,a_{i_{2} j} )\} \le   \max \{ 1-\left|x_{i_{1} } -a_{i_{1} j} \right|\, ,\, 1-\left|x_{i_{2} } -a_{i_{2} j} \right|\, \} \] 
b) If $\min (x_{i_{1} } ,a_{i_{1} j} )\, \ne 0$ and $\, \min (x_{i_{2} } ,a_{i_{2} j} )\ne 0$, then: 
\[\max \{ 1-\left|x_{i_{1} } -a_{i_{1} j} \right|\, ,\, 1-\left|x_{i_{2} } -a_{i_{2} j} \right|\, \} \le S_{D} \{ \min (x_{i_{1} } ,a_{i_{1} j} )\, ,\, \min (x_{i_{2} } ,a_{i_{2} j} )\} \] 
\vskip 0.1in
\noindent\textbf{Proof.} a) Suppose $\min (x_{i_{1} } ,a_{i_{1} j} )\, =\, \min (x_{i_{2} } ,a_{i_{2} j} )=0$. In this case, we have $S_{D} \{ \min (x_{i_{1} } ,a_{i_{1} j} )\, ,\, \min (x_{i_{2} } ,a_{i_{2} j} )\} =0$. Now, by noting that $0\le 1-\left|x_{i} -a_{ij} \right|\le 1$, $\forall i\in I$and$\forall j\in J$, the result is attained. In another case, suppose $\min (x_{i_{1} } ,a_{i_{1} j} )\, =0$ and $\, \min (x_{i_{2} } ,a_{i_{2} j} )\ne 0$. Then, by the definition of drastic sum and Lemma 1, we have:
\[S_{D} \{ \min (x_{i_{1} } ,a_{i_{1} j} )\, ,\, \min (x_{i_{2} } ,a_{i_{2} j} )\} =\, \min (x_{i_{2} } ,a_{i_{2} j} )\le \, 1-\left|x_{i_{2} } -a_{i_{2} j} \right|\, \le \max \{ 1-\left|x_{i_{1} } -a_{i_{1} j} \right|\, ,\, 1-\left|x_{i_{2} } -a_{i_{2} j} \right|\, \} \] 
that implies statement. (b) Suppose $\min (x_{i_{1} } ,a_{i_{1} j} )\, \ne 0$ and $\, \min (x_{i_{2} } ,a_{i_{2} j} )\ne 0$. Then, by the definition of drastic sum, we have $S_{D} \{ \min (x_{i_{1} } ,a_{i_{1} j} )\, ,\, \min (x_{i_{2} } ,a_{i_{2} j} )\} =1$. This fact together with $0\le \max \{ 1-\left|x_{i_{1} } -a_{i_{1} j} \right|\, ,\, 1-\left|x_{i_{2} } -a_{i_{2} j} \right|\, \} \le 1$ prove part b.
\vskip 0.1in
\noindent 

\noindent    Lemmas 1 and 3 lead us to a theorem below by which we will able to introduce and use a new t-conorm instead of maximum operator in (2), for compensating the difference between solution attained from (2) and real solution resulted from \eqref{GrindEQ__4_}.

\noindent 
\vskip 0.1in
\noindent \textbf{Theorem 1.} Consider Problem (2) and suppose for each $j\in J$ there exist at least two $i_{1} ,i_{2} \in I$ such that $\min (x_{i_{1} } ,a_{i_{1} j} )\, \ne 0$and$\, \min (x_{i_{2} } ,a_{i_{2} j} )\ne 0$. Then, $\mathop{\max }\limits_{i=1}^{n} \{ \, \min (x_{i} ,a_{ij} )\, \} \le $$\mathop{\max }\limits_{i=1}^{n} \{ 1-\left|x_{i} -a_{ij} \right|\, \} \le $$\mathop{S_{D} }\limits_{i=1}^{n} \{ \, \min (x_{i} ,a_{ij} )\, \} $.
\vskip 0.1in
\noindent \textbf{Proof.} First inequality is the same which was proved in Lemma 1. Also, second inequality is obtained by repeating Lemma 3.

\noindent 

\noindent   Now, let  $i_{1} ,i_{2} \in I$ such that $\min (x_{i_{1} } ,a_{i_{1} j} )\, \ne 0$and$\, \min (x_{i_{2} } ,a_{i_{2} j} )\ne 0$. By Theorem 1,\\* $\mathop{\max }\limits_{i=i_{1} }^{i_{2} } \{ \, \min (x_{i} ,a_{ij} )\, \} \le $$\mathop{\max }\limits_{i=i_{1} }^{i_{2} } \{ 1-\left|x_{i} -a_{ij} \right|\, \} \le $$\mathop{S_{D} }\limits_{i=i_{1} }^{i_{2} } \{ \, \min (x_{i} ,a_{ij} )\, \} $. \\*Then, for some $\lambda \in [0,1]$, we have: $\lambda \mathop{\max }\limits_{i=i_{1} }^{i_{2} } \{ \, \min (x_{i} ,a_{ij} )\, \} +(1-\lambda )$$\mathop{S_{D} }\limits_{i=i_{1} }^{i_{2} } \{ \, \min (x_{i} ,a_{ij} )\, \} $=$\mathop{\max }\limits_{i=i_{1} }^{i_{2} } \{ 1-\left|x_{i} -a_{ij} \right|\, \} $in which :

\[\lambda =\frac{\mathop{\max }\limits_{i=i_{1} }^{i_{2} } \{ 1-\left|x_{i} -a_{ij} \right|\, \, \} -\mathop{\max }\limits_{i=i_{1} }^{i_{2} } \{ \min (x_{i} ,a_{ij} )\} }{\mathop{S_{D} }\limits_{i=i_{1} }^{i_{2} } \{ \min (x_{i} ,a_{ij} )\} -\mathop{\max }\limits_{i=i_{1} }^{i_{2} } \{ \min (x_{i} ,a_{ij} )\} }  \qquad(*)\] 
   Now, by statements above, we define $S_{a,b} (a,b):\, [0,1]^{2} \to [0,1]$ such that $S_{a,b} (a,b)=$$\lambda \max \{ \min (a,b)\} +(1-\lambda )S_{D} (a,b)$in which$\lambda $ is found by (*). It is easy to verify that $S_{a,b} (a,b)$ is really a t-conorm. Also, by extension of $S_{a,b} (a,b)$for $n$ arguments and usage of it instead of maximum operator in Problem (2), we have:
\[\mathop{S_{M} }\limits_{i=1}^{n} \{ \min (x_{i} ,a_{ij} )\} =\lambda \mathop{\max }\limits_{i=1}^{n} \{ \, \min (x_{i} ,a_{ij} )\, \} +(1-\lambda )\mathop{S_{D} }\limits_{i=1}^{n} \{ \, \min (x_{i} ,a_{ij} )\, \} \] 
in which$M=\{ \min (x_{i} ,a_{ij} )\, ,\, \forall i\in I\, \} $and $\lambda $ is found by equality below:
\[\lambda =\frac{\mathop{\max }\limits_{i=1}^{n} \{ 1-\left|x_{i} -a_{ij} \right|\, \, \} -\mathop{\max }\limits_{i=1}^{n} \{ \min (x_{i} ,a_{ij} )\} }{\mathop{S_{D} }\limits_{i=1}^{n} \{ \min (x_{i} ,a_{ij} )\} -\mathop{\max }\limits_{i=1}^{n} \{ \min (x_{i} ,a_{ij} )\} } \qquad  (**)\] 
   However, if assumption in Theorem 1 hold, statements above show that there exists t-conorm $S_{a,b} (a,b)$ such that $\mathop{S_{M} }\limits_{i=1}^{n} \{ \min (x_{i} ,a_{ij} )\} =$ $\mathop{\max }\limits_{i=1}^{n} \{ 1-\left|x_{i} -a_{ij} \right|\, \} $, $\forall j\in J$. Now, we have proved theorem below.

\noindent 
\vskip 0.1in
\noindent \textbf{Theorem 2.} Consider Problem (2) and suppose for each $j\in J$ there exist at least two $i_{1} ,i_{2} \in I$ such that $\min (x_{i_{1} } ,a_{i_{1} j} )\, \ne 0$and$\, \min (x_{i_{2} } ,a_{i_{2} j} )\ne 0$. Define $S_{a,b} (a,b):\, [0,1]^{2} \to [0,1]$ such that \\* $S_{a,b} (a,b)=$$\lambda \max \{ \min (a,b)\} +$ $(1-\lambda )S_{D} (a,b)$ and 
\[\lambda =\frac{\mathop{\max }\limits_{i=i_{1} }^{i_{2} } \{ 1-\left|x_{i} -a_{ij} \right|\, \, \} -\mathop{\max }\limits_{i=i_{1} }^{i_{2} } \{ \min (x_{i} ,a_{ij} )\} }{\mathop{S_{D} }\limits_{i=i_{1} }^{i_{2} } \{ \min (x_{i} ,a_{ij} )\} -\mathop{\max }\limits_{i=i_{1} }^{i_{2} } \{ \min (x_{i} ,a_{ij} )\} } \] 
Then, $S_{a,b} (a,b)$ is a t-conorm and $\mathop{S_{M} }\limits_{i=1}^{n} \{ \min (x_{i} ,a_{ij} )\} =$$\mathop{\max }\limits_{i=1}^{n} \{ 1-\left|x_{i} -a_{ij} \right|\, \} $,  $\forall j\in J$, in which\\* $M=\{ \min (x_{i} ,a_{ij} )\, ,\, \forall i\in I\, \} $and $\lambda $ is found by
\[\lambda =\frac{\mathop{\max }\limits_{i=1}^{n} \{ 1-\left|x_{i} -a_{ij} \right|\, \, \} -\mathop{\max }\limits_{i=1}^{n} \{ \min (x_{i} ,a_{ij} )\} }{\mathop{S_{D} }\limits_{i=1}^{n} \{ \min (x_{i} ,a_{ij} )\} -\mathop{\max }\limits_{i=1}^{n} \{ \min (x_{i} ,a_{ij} )\} } .\] 
\vskip 0.1in

\noindent \textbf{Theorem 3.} Consider Problem (2) and Suppose for each $j\in J$ there exists at most one $i_{j} \in I$ such that $\min (x_{i_{j} } ,a_{i_{j} j} )\ne 0$. Then, for each t-conorm $S$ and each t-norm T, we have \\*$\mathop{S}\limits_{i=1}^{n} \{ t(x_{i} ,a_{ij} )\} \le $$\mathop{\max }\limits_{i=1}^{n} \{ 1-\left|x_{i} -a_{ij} \right|\, \} $, $\forall j\in J$.

\noindent 
\vskip 0.2in
\noindent \textbf{Proof.} Since  for each $a,b\in [0,1]$ , we have $T(a,b)\le \min (a,b)$ and $S(a,b)\le S_{D} (a,b)$for each t-norm T and each t-conorm $S$, it is sufficient to prove $\mathop{S_{D} }\limits_{i=1}^{n} \{ \min (x_{i} ,a_{ij} )\} <$$\mathop{\max }\limits_{i=1}^{n} \{ 1-\left|x_{i} -a_{ij} \right|\, \} $. At first, suppose that for a fixed $j\in J$, we have $\min (x_{i} ,a_{ij} )\} =0$, $\forall i\in I$. Then, $\mathop{S_{D} }\limits_{i=1}^{n} \{ \min (x_{i} ,a_{ij} )\} =0$ which implies inequality. In another case, suppose for a fixed $j\in J$ there exists exactly one $i_{j} \in I$ such that $\min (x_{i_{j} } ,a_{i_{j} j} )\ne 0$. In this case, $\mathop{S_{D} }\limits_{i=1}^{n} \{ \min (x_{i} ,a_{ij} )\} =\min (x_{i_{j} } ,a_{i_{j} j} )$. Now, by (5) we have:
\[\mathop{S_{D} }\limits_{i=1}^{n} \{ \min (x_{i} ,a_{ij} )\} =\min (x_{i_{j} } ,a_{i_{j} j} )\le \, 1-\left|x_{i_{j} } -a_{i_{j} j} \right|\, \le \mathop{\max }\limits_{i=1}^{n} \{ 1-\left|x_{i} -a_{ij} \right|\, \} \] 
that requires the inequality.

\noindent 
\vskip 0.1in
\noindent \textbf{Corollary 2.} Consider Problem (2) and consider a fixed $j\in J$. If  two following conditions hold, then $\mathop{S}\limits_{i=1}^{n} \{ T(x_{i} ,a_{ij} )\} <$$\mathop{\max }\limits_{i=1}^{n} \{ 1-\left|x_{i} -a_{ij} \right|\, \} $.

\noindent (a) For each $j\in J$ there exists $i\in I$ such that $\max (x_{i} ,a_{ij} )<1$.

\noindent (b) $\min (x_{i} ,a_{ij} )=0$, $\forall i\in I$.
\vskip 0.1in
\noindent \textbf{Proof.} Since $\mathop{S}\limits_{i=1}^{n} \{ T(x_{i} ,a_{ij} )\} \le \mathop{S_{D} }\limits_{i=1}^{n} (\min (x_{i} ,a_{ij} )\} $for each t-norm $T$and each t-conorm $S$, it is sufficient to prove $\mathop{S_{D} }\limits_{i=1}^{n} \{ \min (x_{i} ,a_{ij} )\} <$$\mathop{\max }\limits_{i=1}^{n} \{ 1-\left|x_{i} -a_{ij} \right|\, \} $. By assumption (b), $\mathop{S_{D} }\limits_{i=1}^{n} \{ \min (x_{i} ,a_{ij} )\} =0$. Also, from Theorem 3, we have  $\mathop{S}\limits_{i=1}^{n} \{ t(x_{i} ,a_{ij} )\} \le $$\mathop{\max }\limits_{i=1}^{n} \{ 1-\left|x_{i} -a_{ij} \right|\, \} $. Then, to prove the corollary, it is sufficient to show  $\mathop{S}\limits_{i=1}^{n} \{ t(x_{i} ,a_{ij} )\} \ne $$\mathop{\max }\limits_{i=1}^{n} \{ 1-\left|x_{i} -a_{ij} \right|\, \} $.  By contrary, suppose $\mathop{\max }\limits_{i=1}^{n} \{ 1-\left|x_{i} -a_{ij} \right|\, \} $$=\mathop{S}\limits_{i=1}^{n} \{ t(x_{i} ,a_{ij} )\} =0$. Therefore, $1-\left|x_{i} -a_{ij} \right|=0$, $\forall i\in I$, which implies $\min (x_{i} ,a_{ij} )=$$1-\left|x_{i} -a_{ij} \right|$, $\forall i\in I$. Thus, $\max (x_{i} ,a_{ij} )=1$, $\forall i\in I$, which contradicts assumption (a). Then, the proof is complete.
\vskip 0.1in
\noindent 

\noindent    Briefly, the amount of the accuracy of vector $b$resulted from the direct solution of FRI problem $xo \, A=b$ modeling an applied topic can be categorized as follows:

\noindent 
\vskip 0.1in
\begin{itemize}
\item  If $\max \{ x_{i} ,a_{ij} \} =1$, $\forall i\in I$and $\forall j\in J$, Problem (2) defined by    minimum t-norm yields exact solution(see Lemma 2 part (b)). In this case, Problem (2) is formulated as $xo \, a_{j} \, =\, \mathop{\max }\limits_{i=1}^{n} \{ \, \min (\, x_{i} ,a_{ij} )\, \} \, =\, b_{i} $$,\forall j\in J$.

\item  If for each $j\in J$ there exist at least two $i_{1} ,i_{2} \in I$ such that $\min (x_{i_{1} } ,a_{i_{1} j} )\, \ne 0$and$\, \min (x_{i_{2} } ,a_{i_{2} j} )\ne 0$, then Problem (2) defined by    minimum t-norm and $S_{a,b} $ t-conorm given in Theorem 2, yields exact solution(See Theorems 1 and 2). Here, Problem (2) is formulated as: $xo \, a_{j} \, =\, \mathop{S_{M} }\limits_{i=1}^{n} \{ \, \min (\, x_{i} ,a_{ij} )\, \} \, =\, b_{i} $$,\forall j\in J$.
\vskip 0.1in
\item  However, there are conditions in which we cannot obtain reasonable solutions;
\begin{itemize}
\item  If necessary condition stated in Lemma 2 part (b) does not hold, then Problem (2) defined by minimum t-norm fails to provide a rational solution for vector $b$.

\item  If conditions expressed in Corollary 1 hold, we have: 
\[\mathop{\max }\limits_{i=1}^{n} \{ T(x_{i} ,a_{ij} )\} <\mathop{\max }\limits_{i=1}^{n} \{ 1-\left|x_{i} -a_{ij} \right|\, \} , \forall j\in J\] 

which implies that Problem (2) defined by any t-norm does not yield a justifiable solution $b$.

\item  If conditions in Corollary 2 hold, we have:
\[\mathop{S}\limits_{i=1}^{n} \{ T(x_{i} ,a_{ij} )\} <\mathop{\max }\limits_{i=1}^{n} \{ 1-\left|x_{i} -a_{ij} \right|\, \} , \forall j\in J\] 
\end{itemize}
\end{itemize}
which requires Problem (2) defined by each t-norm and each t-conorm will attain a solution farther than one we expect. Also, in this case, the best solution of   modified Problem (2) is attained by problem below:
\[\mathop{S_{D} }\limits_{i=1}^{n} \{ \min (x_{i} ,a_{ij} )\} =b_{i} , \forall j\in J\] 
Then, the accuracy of each modification of Problem (2) is bounded above.

\noindent 

\noindent   Therefore, Problem (2) may not guarantee that vector $b$ is the same as we expect. Moreover, in a worse case, it may exists an example for which the results of Problem (2) in direct solution is not even a good approximation for expected rational consequences.

\noindent \textbf{}

\noindent \textbf{Conclusion }
\vskip 0.1in
\noindent In this paper, we investigated direct solution of FRE and compared their results with expected real consequences. It was shown that FRE defined by maximum t-conorm and an arbitrary t-norm yields different interpretations. A necessary condition and a sufficient condition were presented that guarantee FRE defined by maximum t-conorm and minimum t-norm attain the same solutions as human mind does. Also, we presented a t-conorm and use it instead of maximum t-conorm in FRE to obtain solutions with highest similarity with real ones. Moreover, we showed that under some conditions, FRE defined by any t-conorm and any t-norm may find a solution which is not reasonable.

\noindent 
\vskip 0.1in
\noindent \textbf{References }
\vskip 0.1in
\noindent [1] Bellman R.E., Zadeh L.A., Decision-making in fuzzy environment, Management Sci. 17(1970).141-164.

\noindent 

\noindent [2] Chen L., Wang P.P., Fuzzy relation equations(I): the general and specialized  solving algorithms, Soft Computing,6(2002)428-435.

\noindent 

\noindent [3] Di Martino F., Sessa S., Digital watermarking in coding/decoding processes with fuzzy relation equations, Soft Computting,10(2006)238-243.

\noindent 

\noindent [4] Di Martino F., Loia V., Sessa S., A method in the compression /decompression of images using fuzzy equations and fuzzy similarities. In: Proceedings of Conference IFSA 2003 (29/6-2/7/2003, Istanbul, Turkey), pp.524-527(2003).

\noindent 

\noindent [5] Di Nola A., Pedrycz W., Sessa S., Some theoretical aspects of fuzzy relation equations describing fuzzy systems, Inform.Sci.34(1984)241-264.

\noindent 

\noindent [6] Di Nola A., Pedrycz W., Sessa S., Fuzzy relation equations and algorithms of inference mechanism in expert systems, in: M.M. Gupta, A. Kandel, W. Bandler and J.B. Kiszka, Eds., Approximate Reasoning in Expert System (Elsevier Science Publishers (North Holland), Amsterdam, 1985)355-367.

\noindent 

\noindent [7] Di Nola A., Pedrycz W., Sessa S., Sanchez E., Fuzzy relation equations theory as a basis of fuzzy modeling: an overview, Fuzzy Sets and Systems 40(1991)415-429. 

\noindent 

\noindent [8] Di Nola A., Russo C., Lukasiewicz transform and its application to compression and reconstruction of digital images, Information Sciences 177(2007)1481-1498.

\noindent 

\noindent [9]  Di Nola A., Sessa S., Pedrycz W., Sanchez E., Fuzzy relation equations and their applications to knowledge engineering, Kluwer, Dordrecht, (1989).

\noindent 

\noindent [10] Dubois D., Prade H., Fuzzy sets and systems: Theory and Applications (Academic Press, New York,1980).

\noindent 

\noindent [11] Fernandez M.J., Gil P., Some specific types of fuzzy relation equations, Information sciences 164(2004)189-195.

\noindent 

\noindent [12]  Guo F.F., Xia Z.Q., An algorithm for solving optimization problems with one linear objective function and finitely many constraints of fuzzy relation inequalities, Fuzzy Optimization and Decision Making, 5 (2006)33-47.

\noindent 

\noindent [13] Han S.C., Li H.X., Wang J.Y., Resolution of finite fuzzy relation Equations based on strong pseudo-t-norms, Applied Mathematics Letters 19(2006)752-757.

\noindent 

\noindent [14] Han S.C., Li H.X., Nots on ``pseudo-t-norms and implication operators on a complete brouwerian lattice'' and ``Pseudo-t-norms and implication operators: direct products and direct product decompositions '' ,Fuzzy Sets and Systems 153(2005)289-294.

\noindent 

\noindent [15] Hirota K., Pedrycz W., Fuzzy relational compression, IEEE Transactions on Systems, Man and Cybernetics, Part B 29\eqref{GrindEQ__3_}(1999)407-415.

\noindent 

\noindent [16] Hirota K., Pedrycz W., Solving Fuzzy relational equations through logical filtering, Fuzzy Sets Syst. 81(1996)355-364.

\noindent 

\noindent [17] Hirota K., Pedrycz W., Specificity shift in solving ,Fuzzy Sets Syst. 106 (1999)211-220.

\noindent 

\noindent [18] Klement E.P., Mesiar R., Pap E.,Triangular norms. Position paper I: basic analytical and algebraic properties,Fuzzy Sets and Systems 143 (2004) 5-26.

\noindent 

\noindent [19] Klement E.P., Mesiar R., Pap E., Triangular norms. Position paper II, Fuzzy Sets and Systems 145(2004)411-438.

\noindent 

\noindent [20] Klement E.P., Mesiar R., Pap E., Triangular norms. Position paper III, Fuzzy Sets and Systems 145(2004)439-454.

\noindent 

\noindent [21] Klir G.J., Yuan B., Fuzzy sets and fuzzy logic: Theory and Applications,  Prentice-Hall,PTR,USA,1995.

\noindent 

\noindent [22] Loia V., Sessa S., Fuzzy relation equations for coding/decoding processes of images and videos, Information sciences 171(2005)145-172.

\noindent 

\noindent [23] Loia V., Pedrycz W., Sessa S., Fuzzy relation calculus in the compression and decompression of fuzzy relations.Inter J Image Graphics,2(2002)1-15

\noindent 

\noindent [24] Markovskii A.V., On the relation between equations with max-product composition and the covering problem, Fuzzy Sets and Systems 153 (2005)261-273.

\noindent 

\noindent [25] Nobuhara H., Pedrycz W., Hirota K., A digital watermaking algorithm using image compression method based on fuzzy relational equations. In: Proceedings of FUZZ-IEEE 2002.IEEE Press, vol 2,pp.1568-1573.

\noindent 

\noindent [26] Nobuhara H., Hirota K., Pedrycz W., Relational image compression: optimizations through the design of fuzzy coders and YUV colour space, Soft Computing 9(2005)471-479.

\noindent 

\noindent [27] Nobuhara H., Hirota K., Di Martino F., Pedrycz W., Sessa S., fuzzy relation equations for compression/decompression processes of colour images  in the RGB and YUV colour spaces(2005).

\noindent 

\noindent [28] Nobuhara H., Pedrycz W., Hirota K., Fast solving method of fuzzy relational equation and its application to lossy image compression. IEEE Trans Fuzzy Sys,8\eqref{GrindEQ__3_}(2000)325-334.

\noindent 

\noindent [29] Nobuhara H., Bede B., Hirota K., On various eigen fuzzy sets and their application to image reconstruction,Information Sciences 176 (2006) 2988-3010.

\noindent 

\noindent [30] Pedrycz W., Fuzzy relational equations with generalized connectives and their applications, Fuzzy Sets and Systems 5(1983)185-201.

\noindent 

\noindent [31] Pedrycz W., On generalized fuzzy relational equations and their applications, Journal of Mathematical Analysis and Applications 107 (1985)520-536.

\noindent 

\noindent [32] Pedrycz W., s-t fuzzy relational equations, Fuzzy Sets and Systems 59 (1993)189-195.

\noindent 

\noindent [33] Pedrycz W., Inverse problem in fuzzy relational equations, Fuzzy Sets  and Systems 36(1990)277-291.

\noindent 

\noindent [34] Pedrycz W., Vasilakos A.V., Modularization of fuzzy relational equations, Soft Computing 6(2002).33-37.

\noindent 

\noindent [35] Peeva K., Kyosev Y., Fuzzy relational calculus, theory, applications and software, World Scientific,2004.

\noindent 

\noindent [36] Perfilieva I., Fuzzy function as an approximate solution to a system of fuzzy relation equations, Fuzzy Sets and Systems 147(2004)363-383. 

\noindent 

\noindent [37] Pedrycz W., An approach to the analysis of fuzzy systems, Int. J. Control 34(1981)403-421.

\noindent    

\noindent [38] Sanchez E., Solution in composite fuzzy relation equations: Application to medical diagnosis in Brouwerian logic, In Fuzzy Automata and Decision Processes,(Edited by M. M. Gupta , G . N . Saridis and B R Games),pp.221-234,North-Holland,New York , (1977).

\noindent 

\noindent [39] Wang S., Fang S.C., Nuttle H.L.W., Solution sets of interval-valued fuzzy relational equations, Fuzzy Optimization and Decision Making 2 \eqref{GrindEQ__1_} (2003)41-60.

\noindent  

\noindent [40] Wang G.J., On the logic foundation of fuzzy reasoning, Inform.Sci. 117 \eqref{GrindEQ__1_} (1999)47-88.

\noindent 

\noindent [41] Wenstop F., Deductive verbal models of organizations, Int. J. Man-Machine Studies 8 (1976).293-311.

\noindent  

\noindent [42] Wu Y.K., Guu S.M., A note on fuzzy relation programming problems with max-strict-t-norm composition, Fuzzy Optimization and Decision Making,3(2004)271-278.

\noindent  

\noindent [43] Xiong Q.Q., Wang X.P., Some properties of sup-min fuzzy relational equations on infinite domains, Fuzzy Sets and Systems 151(2005)393-402.

\end{document}